\documentclass[10pt,onside,reqno]{amsart}

\newcommand{\Cdb}{\mbox{$\mathbb{C}$}}

\newcommand{\Ndb}{\mbox{$\mathbb{N}$}}

\newcommand{\Rdb}{\mbox{$\mathbb{R}$}}

\newcommand{\Zdb}{\mbox{$\mathbb{Z}$}}
\newcommand{\A}{\mbox{${\mathcal A}$}}

\newcommand{\C}{\mbox{${\mathcal C}$}}

\newcommand{\E}{\mbox{${\mathcal E}$}}
\newcommand{\F}{\mbox{${\mathcal F}$}}

\newcommand{\M}{\mbox{${\mathcal M}$}}

\newcommand{\U}{\mbox{${\mathcal U}$}}
\newcommand{\V}{\mbox{${\mathcal V}$}}
\newcommand{\W}{\mbox{${\mathcal W}$}}

\newcommand{\Y}{\mbox{${\mathcal Y}$}}

\theoremstyle{remark}

\theoremstyle{definition}

\begin{document}

\title[Difference and $(\Delta)$ property] {Difference and $(\Delta)$ properties for some new classes}.
\author{Bolis Basit}
       \address {School of Math. Sciences, P.O. Box  28M,
       Monash University, Victoria 3800, Australia}
        \email{ bolis.basit@monash.edu}

\begin{abstract}
{In this paper we study difference and $(\Delta)$ properties for the classes of the form $C_0(J,X)$,  $\frak {g} \U$,  $\U+\frak {g} \V$, where $\U, \V\in \{BUC(J,X), UC(J,X)\}$  and $\frak{g} (t)=e^{it^2}$, $t\in \mathbb{R}$. For functions whose differences belong to  $ \F \in\{C_0(J,X)$, $\frak {g} \U\}$,  we prove a new stronger $(\Delta)$ property (S$\Delta$P):

If $f: J\to X$ and $\Delta_h f \in \F$, $h > 0$, then $f\in C(J,X)$ and $(f-M_hf)\in \F$,  $h > 0$.

\noindent See  (Lemma 2.5, Theorems 3.1, 3.2, Lemma 3.4, Theorem 3.7).
  These results enabled us to prove $(\Delta)$ for   $\U+\frak {g} \V$ even when $\U, \V\in  UC(J,X) $ (Theorem 4.2).  We give a new proof of a theorem of De Bruijn [10] stating: if $J\in \{\Rdb_+, \Rdb\}$,  $\phi\in \Cdb^J$ and $\Delta_s \phi \in C(J,\Cdb)$ for each $s >0$, then  $\phi= G+H$, where $G \in C(J,\Cdb)$ and $H(t+s)=H(t)+H(s)$, $t,s\in J$, for functions $\phi: X^J  $.}
 \end{abstract}
\maketitle

31-10-18

\

\noindent\textbf{0.   Introduction.}

In this paper we revisit the study of functions $f: J\to X$  whose differences belong to a certain class  $\F \subset C(J,X)$. Here $J\in \{\Rdb_+,\Rdb\}$, $X$ is a Banach space. In [10],  De Bruijn  proved that if $\phi\in \Cdb^J$ and $\Delta_s \phi \in C(J,\Cdb)$ for each $s \in J$, then  $\phi= G+H$, where $G \in C(J,\Cdb)$ and $H(t+s)=H(t)+H(s)$, $t,s\in J$. G\"{unzler} extended this result to functions $f\in X^J$ ([8, Proposition 7.3]). Our attempts to prove $(\Delta)$ property for the sum  $\U+\frak {g} \V$, where $\U, \V\in \{BUC(J,X), UC(J,X)\}$ and $\frak{g} (t)=e^{it^2}$, $t\in \mathbb{R}$, led us  to show that $\F\in \{C_0(J,X), \frak{g}BUC(J,X), \frak{g} UC(J,X)\}$ has a stronger $(\Delta)$ property:

If $f: J\to X$ and $\Delta_h f \in \F$, $h > 0$, then $f\in C(J,X)$ and $(f-M_hf)\in \F$,  $h > 0$\quad (S$\Delta$P).

In Theorem 2.4, we give a new proof of  De Bruijn result [10, Theorem 1.3]. In Lemmas 2.5, 3.6 (i), we show that if $\F\in \{C_0 (J,X), \frak{g} UC(J,X)\}$, then $\F$ has (S$\Delta$P). These results seem new even in the case $X=\Cdb$. In Theorem 3.5, we prove that if $f: J\to X$ and $\Delta_s f \in \frak{g} UC(J,X)\}$, $s > 0$, then $f= \frak{g}u+v$, where $u\in BUC(J,X)$ and $v\in \Delta C_0 (J,X)$. In Theorem 4.2, we show that if $f\in L^1_{loc}(J,X)$ and $\Delta_s f \in ( \U +\frak{g} \V)\}$ for each $s > 0$, then $f \in \Delta \U+ \frak{g} BUC (J,X)$. Here $\U, \V\in \{BUC(J,X), UC(J,X)\}$.

 Lemma 2.1, Proposition 2.2 and Lemma 2.3 of  Section 2 are unpublished
   results by the author and H. G\"{u}nzler [6].  In [8],  $(\Delta)$ property for classes  from $\{\frak{g}BUC(J,X), \frak{g} UC(J,X)\}$ is proved  using
  complicated constructions.  Our proofs of (S$\Delta$P) for
   $\F\in \{\frak{g}BUC(J,X),$
$ \frak{g} UC(J,X)\}$ are simple  and this new property allowed us to prove $(\Delta)$ property  even for $ UC(J,X) + \frak{g} UC(J,X)$ (see  [8, 9. Open questions, no 8, p. 30]).  I would like to thank Hans G\"{u}nzler for his numerous  corrections to this paper.

\

\noindent\textbf{1.   Notation, definitions.}

In the following   $\mathbb{R}_+ = [0,\infty)$, $J\in \{\mathbb{R}_+, \mathbb{R}\}$,  $\Ndb=\{1, 2, \cdots\}$ and $(X, ||\cdot||)$ is a real or a complex Banach space.  We denote by
\begin{equation}
\label{metric}
X_c,  C(J,X), UC(J,X) , BC(J,X), BUC(J,X), C_0(J,X),
\Delta
\A (J,X),
\quad  J\in\{\Rdb_+,\Rdb\}
\end{equation}
\noindent  the spaces of constant,  continuous,  uniformly  continuous,  bounded continuous,  bounded uniformly  continuous, continuous vanishing at infinity functions and
\begin{equation}
\label{metric}\Delta  \A (J,X):=\{\phi\in L^1_{loc}(J,X): \Delta_h \phi \in \A (J,X), h >0 \}.\end{equation}
\noindent Endow $C (J,X)$ with the infinite sequence of  semi-norms
\begin{equation}
||f||_n=||f| [-n,n]\cap J||_{\infty} \text{\,\,with \,\,}n\in \Ndb \text{\,\, and\,\,}  f\in C (J,X).
\end{equation}
\noindent Then   $C (J,X)$ is a metriziable
 complete locally convex linear topological  space [13, p. 31].

\

\noindent If  $\F\in \{BUC(J,X), BC(J,X)$, $e^{i\,t^2} BUC(J,X)$, $UC(J,X),e^{i\,t^2}UC(J,X)\}$, then $\F$ is a  linear subspace of $C (J,X)$ but $\F$ is not closed  in $C (J,X)$ endowed with topology of uniform convergence on compact intervals defined by semi-norms $||f||_n$. But  $\F$ is a complete locally convex linear topological    space endowed with the topology of uniform convergence on $J$ defined by the  neighbourhood system with base
\begin{equation}
 \{V_{1/n}(f): f\in \F, n\in \Ndb\},\text{\, where \,} V_{1/n}(f)= \{h\in \F: ||h-f||_{\infty}\le 1/n\}.
\end{equation}
\noindent Moreover, $BUC(J,X), BC(J,X), e^{i\,t^2} BUC(J,X)$
 are Banach subspaces
 and  $ C_0 (J,X)$ is a Banach subspace of $BUC(J,X)$.

\noindent Set
\begin{equation}
\label{metric}
\frak{g}(t) = e^{i\,t^2} ([9, p. 202, \text{\,problem\,}  8.18]),\quad  \gamma_{2s}(t)= e^{2i\, st}, \quad E(\frak {g})= \frac{X}{\frak {g}}, t\in \Rdb.
\end{equation}
\noindent Obviously $E(\frak {g})$ is a closed subspace of $BC (\mathbb{R},X)$. Define
\begin{multline}
\label{metric}
(T\phi) (s,t)=\frak{g}(s)\gamma_{2s}(t)\phi(s+t)- \phi(t)=(1/\frak{g}(t))\Delta_s (\frak{g}(t)\phi(t)),\,\,
\phi\in  X^J \text{\, and\,} s,t\in J,\\
(\underline{T}\phi) (s,t)=\overline{\frak{g}(s)\gamma_{2s}(t)}\phi(s+t)- \phi(t)=\frak{g}(t)\Delta_s (\phi(t)/\frak{g}(t)),\,\,
\phi\in  X^J \text{\,\, and\,\,} s,t\in J.
\end{multline}
\noindent Then $T:  C(\mathbb{R},X)\to  C(\mathbb{R}^2,X)$ is a linear continuous operator, $ T E(\frak {g})=\{0\}$   and  $T  C(\mathbb{R},X) \subset Z_0$, where

\

 $Z_0= :\{\psi\in C(\mathbb{R}^2,X): \psi (0,t)=0 \text{\,\, for\,\,} t\in\mathbb{R}\}$, so
 $Z_0$  is a closed subspace of $C(\mathbb{R}^2,X)$.

\

\noindent If  $\A, U  \subset  X^J$, we say that
$\A$ $satisfies$  $(L_U)$  (see [3, p. 677]) if
        $F  \in  U$ with all differences  $\Delta_h F \in \A$  for  $h>0$ imply $F\in\A$.
        $\A\subset L^1_{loc} (J,X)$ (see [1, p. 13])  $satisfies$ $(P_U)$ if  $f \in \A$  with indefinite integral  $Pf \in U$  implies  $Pf \in \A$.
        $(P_b)$ means $(P_{C_b})$, similarly for $(L_b)$.

\

\noindent We say that $\A  \subset   L^1_{loc}(J,X)$  $satisfies$
$(\Delta)$  (see [3, p. 677]) if for any $f \in L^1_{loc}(J,X)$ for which all
differences $\Delta_s f \in  \A$ for $0< s \in \Rdb$  one has $f -
M_h f  \in \A$ for all $h>0$, where
\begin{equation}
\label{metric}
M_h f (t)=(1/h) \int_{0}^h f(t+s)\, ds.
\end{equation}

\noindent\textbf{2.  Preliminaries.}
In this section  we prove some needed  results which might be of independent interest. We study the linear subspaces $UC(J,X)$, $\frak{g}\cdot UC(J,X)\subset  C(J,X)$. We endow $UC(J,X)$, $\frak{g}\cdot UC(J,X)$ with the topology of uniform convergence on $J$ given by the base neighbourhood system (4). Then  $UC(J,X)$, $\frak{g}\cdot UC(J,X)$ are complete locally convex linear spaces.

\

Let $\mathcal {X}_1$ and $\mathcal {X}_2$ be two topological linear spaces over the same field of scalars. Let $\mathcal {X}_1\oplus \mathcal {X}_2$ be the direct sum of the linear spaces $\mathcal {X}_1$ and $\mathcal {X}_2$ in the sense of [11, Section I.11, p. 37] with the product topology of [11, Section I.8, p. 31]. Then $\mathcal {X}_1\oplus \mathcal {X}_2$ is readily seen to be a topological linear space.
If  $\mathcal {X}_1$ and $\mathcal {X}_2$ are Banach- (or F-) spaces, then  $\mathcal {X}_1\oplus \mathcal {X}_2$ is a Banach- (or an F-) space under either of the norms
\begin{multline}
\label{metric} ||x+y||=\max\{||x||,||y||\}, \quad ||x+y||= (||x||^p+||y||^p)^{1/p}, x\in \mathcal {X}_1,  y\in \mathcal {X}_2,\,\,  1\le p < \infty,\\
\text{\,\, and\, these\, norms\, are\, equivalent\, to \,the \, product topology}.\qquad\qquad\qquad\qquad\qquad
\end{multline}

\noindent\textbf{Lemma 2.1}. (i) If $u, v \in UC(J,X)$, then

 $ {\lim}_{T\to \infty}
 \sup\{||u(t)||: |t| \ge T\}$,\,
              $ {\lim}_{T\to \infty}
 \sup\{||v(t)||: |t| \ge T\}\le ||u +\frak{g}^k v
||_{\infty}$, $k=1, 2$.

\noindent (ii) If $u, v \in UC(J,X)$ with  $||u+g^k v||_{\infty}
< \infty$, $k=1, 2$, then $u, v$ are bounded.
           So,

             $ ((UC(J,X) +\frak{g}^k UC(J,X)) \cap BC (J,X) =  BUC (J,X) +\frak{g}^k BUC (J,X)$, $k =1, 2$.

\noindent (iii) If $u, v \in UC(J,X)$, there exists $y\in C_0 (J,X)$ such that $u+\frak{g} v=u-y +\frak{g}(v+y/\frak{g})$ satisfies  $||u-y||_{\infty}\le  ||u+\frak{g} v||_{\infty}$.

\noindent (iv) $ (UC (J,X) \cap  \frak{g}^k UC (J,X))  =BUC (J,X) \cap
\frak{g}^r BUC (J,X)
= C_0 (J,X)$, $k, r =1, 2$.

\noindent (v) Let $ \U, \V\in \{BUC(J,X), UC(J,X)\}$ and $ \W=\U+\frak {g} \V$. Then $\W$ is a complete locally convex linear space endowed with the topology of uniform convergence on $J$.
\

\noindent\textbf{Proof.} (i) It is enough to prove the case $X=\Cdb$. We have
$||u+\frak{g}^k v||_{\infty}\ge|u(t)+\frak{g}^k(t)v(t)|$ for all $t\in J$, $k=1, 2$.
Since $\frak{g}^k (t)$ assumes all complex values $\{z\in\Cdb: |z|=1\}$, $k=1, 2$, it follows $|u(t')+\frak{g}^k (t')v(t')|=|u(t')|+|v(t')|$
for some $t'\in [t,t+\delta (t)]$ for $0<t\in\Rdb$, where $\delta (t) =\sqrt {t^2+ 2\pi}-t$.
 Since $u,v \in UC (J,X)$,  for each $\varepsilon >0$, there exists $\delta (\varepsilon) > 0$ such that $|u(s')-u(s'')|\le \varepsilon/2$,  $|v(s')-v(s'')|\le \varepsilon/2$ if $|s'-s''|\le \delta (\varepsilon) $. Choose  $T>0$ such that $\delta (t) <\delta (\varepsilon) $ if $|t| \ge T$. It follows $||u+\frak{g}^k v||_{\infty}\ge |u(t)|+|v(t)|-\varepsilon$ for all $|t|\ge T$. Part (i) follows.

 (ii) By part (i), for each $\varepsilon > 0$, there is $T > 0$ such that $||u(t)||+||v(t)||-\varepsilon \le ||u+g^k v||_{\infty} < \infty$ for all $|t|\ge T$, $k=1, 2$. Since $u, v$ are bounded if $|t| \le T$, part (ii) follows.

(iii) Choose $y\in C_0 (J,X)$ such that $||u-y||_{\infty}=  {\lim}_{T\to \infty}
 \sup\{||u(t)||: |t| \ge T\}$.  By part (i), we get  $||u-y||_{\infty}\le ||u-y+\frak{g} (v+y/\frak{g})||_{\infty}= ||u+\frak{g} v||_{\infty}$.

 (iv)   By part (ii) it is enough to prove  $BUC (J,X) \cap
\frak{g} BUC (J,X)= C_0 (J,X)$. Let $\psi= \frak{g}\phi\in \frak{g}\cdot BUC(J,X)
\cap BUC(J,X)$.  Assuming $\psi\not\in C_{0}(J,X)$, there is
a sequence $(t_n)\subset J$  such that $|t_n|\to \infty$, and $||\psi
(t_n)||
> \varepsilon_0 > 0$  for all $n\in \mathbb{N}$.  By triangular inequality, one
gets $||\frac {\psi(t_n+h)}{\frak{g}(t_n+h)}- \frac
{\psi(t_n)}{\frak{g}(t_n)}||\ge |||\frac
{\psi(t_n+h)}{\frak{g}(t_n+h)}- \frac {\psi(t_n
+h)}{\frak{g}(t_n)}|| - ||\frac
{\psi(t_n+h)-\psi(t_n)}{\frak{g}(t_n)}|||$ for all $n\in \mathbb{N}$ and
$h\in J$. Using uniform continuity of $\psi$ and rapid
oscillation of
 $\frak{g}$, there is a sequence $(h_n)$ such
that $h_n\to 0$ as $n\to\infty$  and $||\frac
{\psi(t_n+h_n)}{\frak{g}(t_n+h_n)}- \frac
{\psi(t_n)}{\frak{g}(t_n)}||\ge \varepsilon_0/4$  for all $n\ge m$ for some
$m\in \mathbb{N}$. This implies $\phi\not \in BUC(J,X)$, a
contradiction which proves part (iii).

(v)  By part (iv), $\U \cap \frak{g} \V= C_0 (J,X)$.  By part (i), if  $u, v \in UC(J,X)$ and $||u+\frak{g} v||_{\infty} < \infty$, then  $u,v \in BUC(J,X)$. By part (iii), if $(u_n+\frak{g} v_n)$ is a Cauchy sequence, then $(u_n)$ can be replaced by $(u_n-y_n)$ with $(y_n) \subset C_0 (J,X)$ so that $u'_n=u_n-y_n$,  $v'_n=v_n+y_n/\frak{g}$ satisfy $|| u'_n-u'_m||_{\infty} \le || (u'_n-u'_m)+\frak{g} (v'_n-v'_m)||_{\infty}= || (u_n-u_m)+\frak{g} (v_n-v_m)||_{\infty}$ for all  $n, m\in \Ndb$. It follows  $(u_n')$,   $(\frak{g}v'_n)$ are Cauchy sequences. Since $\U,  \frak{g} \V$ are locally convex complete linear spaces,  $(u_n+\frak{g} v_n)$ converges to an element from $\U + \frak{g} \V$.
 $\square$

\

\noindent\textbf{Proposition 2.2}. Let $ \U, \V\in \{BUC(J,X), UC(J,X)\}$ and $ \W=\U+\frak {g} \V$. Then
 $\W/C_{0}(J,X)= $
 $ \U/C_{0}(J,X)\oplus (\frak{g}\cdot \V)/C_{0}(J,X)$ is a  locally convex complete linear space.

\

\noindent\textbf{Proof.} As mentioned above  $\U,\frak{g} \V$ are  locally convex complete linear spaces and by  Lemma 2.1 (iv), $\W$  is a locally convex complete linear space. Since  $\C_0(J,X)$ is a closed subspace of $\U,\frak{g} \V, \W$, the quotient spaces   $\W/C_{0}(J,X) $,
 $ \U/C_{0}(J,X)$, $ (\frak{g}\cdot \V)/C_{0}(J,X)$  are  locally convex complete linear spaces ([12, (2), p. 207]). By Lemma 2.1 (iii),  $ \U/C_{0}(J,X)\cap (\frak{g}\cdot \V)/C_{0}(J,X)= \{0\}$. Since  $ \U/C_{0}(J,X), (\frak{g}\cdot \V)/C_{0}(J,X)\subset \W/C_{0}(J,X) $, we conclude  $\W/C_{0}(J,X)= $
 $ \U/C_{0}(J,X)\oplus (\frak{g}\cdot \V)/C_{0}(J,X)$. $\square$

\

\noindent\textbf{Lemma 2.3}. (i) $\frak{g}\cdot BUC(J,X)\subset \M  C_{0}(J,X)\subset\E_0 (J,X)$.

(ii) $M_h (\frak{g}\gamma_{\lambda}\phi) \in \frak{g} BUC(J,X)$ for each $\gamma_{\lambda}\in \widehat {\Rdb}$, $\phi \in UC(J,X)$.

(iii) $ \frak{g}UC (J,X) \subset   \M(\frak{g}UC (J,X)$.

\

\noindent\textbf{Proof.} (i) Let $\psi= \frak{g}\phi\in \frak{g}\cdot
BUC(J,X)$. Then $M_h \psi\in BUC(J,X)$ for each $h>0$. We
have  $hM_h \psi (t)=\frak{g}\int_0^h
\frak{g}(s)\gamma_{2s}(t)\phi(t+s)\,ds=\frak{g}F$ with $F\in
BUC(J,X)$. It follows $hM_h \psi\in BUC(J,X)\cap
\frak{g}\cdot BUC(J,X) =C_0 (J,X)$ by Lemma
 2.1  (iv). This proves  $\frak{g}\cdot BUC\subset \M  C_{0}$.

Now, let $f\in \M C_0 (J,X)$.  Then $M_1  f\in C_0 (J,X)$,
$M_n f= (1/n) \sum_{k=0}^{n-1} (M_1 f)_{k}$. Direct calculations show that $M_n f\to 0$ as $n\to \infty$. This implies $f\in \E_0 (J,X)$.

(ii) Let $\phi \in UC (J,X)$. We prove the case $\gamma_{\lambda}=1$ for the other cases are similar.   We have

 $M_h (\frak{g}\phi) (t)= I_1 +I_2$, where

$ I_1= (1/h) \frak{g}(t) \int_0^h \frak{g}(s)\gamma_{2s} (t)(\phi (t+s)-\phi (t)) ds $, $I_2 =(1/h) \frak{g}(t) \phi(t) \int_0^h \frak{g}(s)\gamma_{2s} (t) ds$.

\noindent   Since $\phi \in UC$,   $(\phi_s-\phi) \in BUC$ uniformly on $[0,h]$. It follows  $I_1 \in g  BUC$.
 Integrating  $\int_0^h \frak{g}(s)\gamma_{2s} (t) ds$ by parts if  $|t| \ge 1$,

       $(h/\frak{g}(t))  I_2 =
        (\phi(t)/{2it})[( \frak{g}(s)\gamma_{2h}(t) -1)-2i \int_0^h s \, \frak{g}(s)\gamma_{2s} (t) ds]$.

\noindent Direct calculations show that $||\phi (t)|| \le ||\phi(0)||+ c |t| $, $t\in J$  for some constant $c > 0$. It follows $\phi(t)/{2it}$,  $(\frak{g}(s)\gamma_{2h}(t) -1)$, $ -2i \int_0^h s \, \frak{g}(s)\gamma_{2s} (t) ds$ are uniformly continuous bounded functions if  $t\in  J$, $|t| \ge 1$. Since  $I_2$ is uniformly continuous bounded on $[-1,1]\cap J$, we conclude $I_2 \in \frak{g} BUC(J,X)$.

(iii) Since  $\frak{g} BUC(J,X) \subset \frak{g} UC(J,X)$, part (iii) follows by (ii).
$\square$

\

\noindent\textbf{Theorem 2.4 (DeBruijn [10, Theorem  1.3])}. If $ \phi \in \Cdb^J$ and $\Delta_s \phi\in C(J,\Cdb)$ for each $s > 0$, then $\phi (t) = G(t)+ H(t)$, where $G\in C(J,\Cdb)$ and $H (t+s)=H(t) +H(s)$ for each $t, s \in J$.

\

\noindent\textbf{Proof}.
Assuming $\phi$ is bounded  on a neighbourhood of some point $t_0\in J$, then by [2, Theorem 4.1], $\phi$ is continuous at $t_0$. It follows $\phi$ is continuous on $J$. It means that either $\phi$ is continuous or $\phi$ is unbounded on each neighbourhood of every point $t\in J$. Write $\Cdb^J/C(J,\Cdb)= \mathcal {D}+ \mathcal {A}$, where $\mathcal {A}=\{H\in \Cdb^J/C(J,\Cdb): H (t+s)=H(t)+H(s),  t,s \in J\}$. We note that if $\phi+ C(J,\Cdb)=\psi\not\in \mathcal {A}$ and $\psi$ is unbounded on each neighbourhood of every point $t\in J$, then $\psi$ is not additive, so $ \psi_s -\psi \not =  \psi (s)$ for some $0\not = s\in J$. Hence $\psi_s -\psi $ is not continuous. This contradicts the assumptions. Hence $\psi\in \mathcal {A}$. It follows that $\phi\in C(J,\Cdb)+ H$ for some $H$ with $H(t+s)=H(t)+ H(s)$ for each $t, s \in J$.
$\square$

\

\noindent For the proof  of DeBruijn see
 [10, Theorem 1.3] and also the remarks in the proof of [8, Proposition 7.2, p. 26).

\

Since  the spaces $\Cdb, E_{\frak{g}}= (1/\frak{g})\Cdb$ and the space of linear continuous functions  $L_c = \{ l\in C (J,\Cdb): l= ct$ for some $c\in \Cdb\}$ are one dimensional subspaces, there are   subspaces  \,\, \,\,\, $C_{\Cdb}(J,\Cdb), C_{\frak{g}} (J,\Cdb)$,
\newline $C_{L_c} (J,\Cdb) \subset C (J,\Cdb)$ of codimension $1$ such that

\begin{equation}
C (J,\Cdb) = C_{\Cdb}(J,\Cdb)\oplus\Cdb = C_{\frak{g}} (J,\Cdb)\oplus E_{ \frak{g}} =  C_{L_c} (J,\Cdb)\oplus L_c.
\end{equation}

\noindent We note that $U C (J,\Cdb)\cap E_{\frak{g}}=\{0\}$ and $\frak{g} U C (J,\Cdb)\cap \Cdb=\{0\}$, so

 $U C (J,\Cdb) \subset  C_{\frak{g}} (J,\Cdb)$\quad  and \quad
 $\frak{g}U C (J,\Cdb)\subset C_{\Cdb}(J,\Cdb) $.

\

\noindent\textbf{Lemma 2.5}.  (i) If $ \phi\in X^J$ and  $\Delta_h \phi \in C_0 (J,X)$ for each $ h >0$, then $\phi\in UC(J,X)$.

(ii)  If $ \phi\in X^J$ and  $\Delta_h (\frak{g} \phi) \in\frak{g}  UC (J,X)$ for each $ h >0$, then $\phi\in C(J,X)$.

\

\noindent\textbf{Proof}.
(i) First, we prove the case $X=\Cdb$. By  Theorem 2.4 (see
 also [8, p. 26, after (7.2)]), we have $\phi (t) = G(t)+ H(t)$, where $G\in C(J,\Cdb)$ and $H (t+s)=H(t) +H(s)$ for each $t, s \in J$. It follows
$\Delta_h \phi (t)= \Delta_h G (t) + H(h)$. This implies $\Delta_h G \in BUC(J,\Cdb)$ for each $h \in J$. By [2, Corollary 5.5], $G\in UC(J,\Cdb)$.
Hence $\lim _{h\to 0}\Delta _h G=0$ uniformly and so  $\lim_{h\to 0} H(h)= -\lim_{h\to 0}(\lim_{|t|\to \infty}\Delta_ h G (t))
=0$ (see [11, p. 28, 6. Lemma]). It follows $H\in C(J,\Cdb)$ and hence $H\in UC(J,\Cdb)$. This gives $\phi \in UC(J,\Cdb)$.

The case $X$ is a Banach space. By the above $x^* \phi \in UC(J,\Cdb)$ for each $x^* \in X^*$.   It follows $\phi$ is weakly continuous and hence it is weakly locally bounded on $J$ and so locally bounded on $J$  by [11, p. 66, 20. Theorem]. Again, by [2, Corollary 5.5], $\phi \in UC(J,X)$.

(ii) First, we prove the case $X=\Cdb$. As in part (i),  $(\frak{g}\phi) (t) = G(t)+ H(t)$, where $G\in C(J,\Cdb)$ and $H (t+s)=H(t) +H(s)$ for each $t, s \in J$. We can assume that   $G\in C_{L_c}(J,\Cdb) $ defined by (9?).  We have
$\Delta_h(\frak{g} \phi)= (\Delta_h G  + H(h))\in \frak{g} UC (J,X)$ for each $ h >0$. Hence $ \Delta_h G  = \frak{g} u(h,\cdot)- H(h) $,  where $u(h,\cdot)\in  UC (J,X)$  for each $h > 0$.  If $H(h)\not = 0$ for some $h >0$, then  $\Delta_h G\not \in C_{\Cdb} (J,\Cdb)$ by the note below (9?). It follows $G \not\in C_{L_c} (J,\Cdb)$. This is a contradiction which proves $H(h) = 0$ for each $h > 0$. It follows $H=0$ and $\phi = (1/\frak{g})G\in C(J,\Cdb)$.

The case $X$ is a Banach space. By the above $x^* \phi \in  C(J,\Cdb)$ for each $x^* \in X^*$.   It follows $\phi$ is weakly continuous.   So as in part (i), $\phi$ is locally bounded on $J$.  By [4, Proposition 1.3, p. 1009], $\phi \in C(J,X)$.
$\square$

\

 \noindent {\textbf {Remark} 2.6}. The condition $\Delta_h (\frak{g} \phi) \in\frak{g}  UC (J,X)$ for each $ h >0$ in Lemma 2.5 (ii) can not be replaced by  $\Delta_h (\frak{g} \phi) \in  BUC (J,X)$ for each $ h >0$. Indeed, $\phi= (1/\frak{g}) H$, where $H$ is nonmeasurable additive function on $J$,  satisfies the condition $\Delta_s \frak{g}\phi= \Delta_s H= H(s)\in BUC(J,X)$ for each $s > 0 $, but $\phi$ is not continuous.

\

  \noindent {\textbf {Remark} 2.7}. H. G\"{u}nzler gave the following proof of Lemma 2.5 (i):
   
   By  [8, Proposition 7.3], there are  $G\in C(\Rdb,X)$ and additive $H : \Rdb \to X$ with   $\phi  = G + H$  on $\Rdb$; $\Delta_h f  \in  C_0 (\Rdb,X)$ implies $\Delta_h G  \in  Lim : = Lim(\Rdb,X) : =$
$\{F \in  C(\Rdb,X) : \lim_{|t| \to \infty} F(t)$ exists $\in X \}$   (so $\lim_{t \to \infty}  F (t)  =  \lim_{t \to -\infty} F (t)$).

\noindent $Lim$ is a linear invariant closed subspace of $BUC$, $Lim$ satisfies $(\Delta)$  ([8, Proposition 4.1]); [3, Proposition 3.3] gives with $X \subset Lim $:  $G (t)= p(t) + \int _0^t q(s)\, ds$ with some $p, q \in Lim$.

Then $ \Delta_h \phi  =  \Delta_h G +  H(h) = \Delta_h p + h M_h q + H(h)$, $t \to \infty$ gives  $H(h) = h\, lim_{t\to \infty} q (t)$  on $ \Rdb^+$, so $H$ is continuous on $\Rdb^+$ and then on $\Rdb$.

So $ f \in L^1_{loc} (\Rdb,X)$, $f \in$  $\Delta C_o =  C_o + X + P C_o $ by [8, (4.19)].

\

\noindent\textbf{3. Difference and $(\Delta)$ properties for $\frak{g}U$, where $U\in \{BUC(J,X), UC(J,X)\}$.}

\

 Let $(X,||\cdot||_{X})$ be a complex Banach  space  and that $X_r $ is the same space considered as a real linear space  (see [12, p. 273]). Consider  the complex linear space $E=X_r\oplus i X_r$ with multiplication of an element $(y_1+ iy_2)\in E$ by a complex number $a+ib$   defined  by

   $(a+ib)(y_1+ iy_2)=(ay_1-by_2)+i (by_1+ ay_2) $ (see [12, p. 208], [7, p. 134]).

\noindent  Then $E$ is  the direct sum $X_r\oplus i X_r$ which is  isomorphic to the product $X_r\times X_r$ identifying $X_r$ with $(X_r,0)$ and $i X_r$ with $(0,X_r)$.

\noindent Endow $E$ with the  norm $||\cdot||$ defined by (8?) with $X_r =\mathcal {X}_1$ and $iX_r =\mathcal {X}_2$.

Set  $||x_1+i x_2||_r= \sup_{f\in X_r^*, ||f||\le 1} |f(x_1)+i f(x_2)|$, $(x_1+i x_2)\in E$.

\noindent Note that each element  $(x_1+i x_2) \in E$ is also an element of $X=  X_r$.
 By [12, (8), p. 208], [11, 5\, Theorem,  p. 58] (see also  [7, Lemma 7, (72),  p. 134]): With $p=1$ of (8?) and $x_1 \in \mathcal {X}_1$,$i x_1 \in \mathcal {X}_2$,

 $||x_1+i x_2||_r= ||x_1+i x_2||_{X}\le ||x_1||_{X}+||x_2||_{X}= ||x_1+i x_2||$,

$(E,||\cdot||_{r})$ is equivalent to the Banach space  $(E,||\cdot||)$.

\

 To avoid repetitions, in the following

\begin{equation}
\label{metric}
 X \text { \, is\, a\, real\,Banach\, space, \,} E=X\oplus i X \text{\, is\, the\, complex\, space \, as\,  above}.\qquad\,
\end{equation}

\noindent Denote by $BE^J= B(J,E)$ the set of all bounded functions $f: J\to E$.

Using (8?) and (10) with $X= \mathcal{X}_1$,  $iX= \mathcal{X}_2$, it is easy to show

\begin{multline}
\label{metric}
\text {(i)\,}\qquad B(J,E)=B(J,X)\oplus i B(J,X),\qquad \, E^J= X^J\oplus i X^J, \\
\text {(ii)\, if \,} f,k \in X^J, \phi =  (f+i k)\in\Y\in\{BUC(J,E), UC(J,E), C(J,E)\}, \text{\, then \,} f,k\in \Y(J,X).\\
\text{(iii)\,} \,\, BUC(J,E)=BUC(J,X)\oplus iBUC(J,X), \quad\, UC(J,E)= UC(J,X)\oplus i UC(J,X),\\ \text{(iv)}\qquad\qquad\qquad \quad C(J,E)= C(J,X)\oplus i C(J,X).\qquad\qquad\qquad\qquad\qquad\qquad\qquad
\end{multline}

\noindent Indeed,  case $\Y = UC(J,E)$: if $||\phi (t+h)-\phi (t)||= ||(f(t+h)- f(t))+i (k(t+h)-k(t))|| \to 0$ as $|h|\to 0$ uniformly in $t\in J$, then by  (8?), $||(f(t+h)- f(t))||\to 0$,  $||(k(t+h)-k(t))||\to 0$ as $|h|\to 0$ uniformly in $t\in J$.

\noindent  Cases $\Y= BUC(J,E)$ and  $\Y= C(J,E)$  are similar.

\

\noindent{\bf {Lemma 3.1}}. (i) If $\phi \in \Rdb^J$
and
 $(T\phi)(s,\cdot)$ $\in UC (J,\Cdb)$ for each $s > 0
$, then $\phi\in BUC (J,\Rdb)$.

(ii) Let  $X$, $E$ be defined by  (10). If  $\phi \in X^J$
and
 $(T\phi)(s,\cdot)$ $\in UC (J,E)$, $s > 0$, then $\phi\in BUC (J,X)$.
\

 \noindent{\bf{Proof}}. By Lemma 2.5(ii), $\phi\in C(J,X)$.

 (i)  Since $(T\phi)(s,\cdot)$ $\in UC (J,\Cdb)$ for each $s > 0$, $\Delta_{\eta}(T\phi)(s,\cdot)=((T\phi)(s,\cdot+\eta)-(T\phi)(s,\cdot))$ $\in BUC (J,\Cdb)$ for each $s, \eta >0$.
 We have
\begin{multline}
\label{metric}
\,\, \Delta_{\eta}(T\phi)(s,t)=V(s,t)+\gamma_{2s} (t) (\gamma_{2s}(\eta)-1) \phi (t+s), \text{\, where\,}\\
\quad\quad\,\, V(s,t)=\frak{g}(s)\gamma_{2s} (t+\eta) \, (\phi (t+s+\eta)-\phi (t+s))-
(\phi (t+\eta)-\phi (t)).\qquad\qquad\
\end{multline}
\noindent Let $\eta  >0$. Choose $s> 0$ such that $\gamma_{2s}(\eta)=1$. Then
\begin{equation}
\label{metric} \Delta_{\eta}(T\phi)(s,t)=\frak{g}(s)\gamma_{2s} (t+\eta)
 (\phi (t+s+\eta)-\phi (t+s))-(\phi (t+\eta)-\phi (t))=I_1 +i I_2.
\end{equation}
\noindent Since $\Delta_{\eta}(T\phi)(s,\cdot)\in BUC (J,\Cdb) $, by (13) and (11) (iii), $I_1, I_2 \in BUC(J, \Rdb)$, where

 $I_1=(\cos (s^2 +2s(\cdot+\eta))\,\,(\phi_{s+\eta} -\phi_s)
-(\phi_{\eta} -\phi))$,\,\,
 $ I_2=\sin  (s^2 +2s(\cdot+\eta))\,\,(\phi_{s+\eta} -\phi_s) $.

\noindent Let $\varepsilon >0$,  $K(\varepsilon,s)= \{t\in J : |\sin  (s^2 +2s(t+\eta))| < \varepsilon\}$. Since $I_2 \in BUC(J, \Rdb)$,
 $(\phi_{s+\eta} -\phi_s $ is bounded and uniformly continuous on $J\setminus K(\varepsilon,s)$ for each  $\varepsilon > 0$. It follows $ \cos (s^2 +2s(\cdot+\eta) (\phi_{s+\eta} -\phi_s) $ is bounded and uniformly continuous on $J\setminus K(\varepsilon,s)$ for each  $\varepsilon > 0$.
  This and  $I_1 \in BUC(J, \Rdb)$  give
\begin{equation}
\label{metric}
 \phi_ {\eta}-\phi\text{\,\, is \, bounded\, and\, uniformly \, continuous\, on \,\,} J\setminus K(\varepsilon,s), \,\,   \varepsilon > 0.
 \end{equation}
\noindent Take $\eta > 0$  and $2s\eta =2\pi$. Then $s=\pi/\eta$. Simple calculations show that one can
  choose $\varepsilon > 0$ such that $(K(\varepsilon,\pi/\eta)-\pi/\eta) \subset J\setminus K(\varepsilon,\pi/\eta)$. This implies
\begin{equation}
\label{metric}
 (\phi_{\eta}-\phi)\text{\,\, is \, bounded\, and\, uniformly \, continuous\, on \,\,}
    K(\varepsilon).
   \end{equation}
 \noindent Combining (14) and (15) together,
 $(\phi (\cdot +\eta)-\phi (\cdot))$ is bounded  and uniformly continuous on $J$. Since $\eta > 0$ is arbitrary, it follows
 $V(s,t) $  is bounded and uniformly continuous on $J$ for each $ s  >0$.
By (12), we conclude   $\gamma_{2s}  (\gamma_{2s}(\eta)-1)(\phi (\cdot+s)=\Delta_{\eta}(T\phi)(s,t) -V(s,t) $  is bounded and uniformly continuous on $J$ for all $s, \eta >0$. This implies $\phi_s $  is bounded and uniformly continuous on $J$  for each $s >0$ and hence $\phi\in BUC (J,\Rdb)$.

(ii) Follows with the same argument as in part (i).
$\square$

 \

\noindent{\bf {Theorem 3.1}
}. (i) Let $\phi: J\to E $ with $E$ defined by (10) and $(T\phi)(s,\cdot)$, $(\underline{T}\phi)(s,\cdot)$ $\in UC(J,E)$ for each $s >0$. Then $\phi\in BUC(J,E)$.

(ii) Let $\phi: J\to E
$  and $(T\phi)(s,
 \cdot)\in C_0 (J,E)$  for each $s >0
$. Then $\phi\in  (1/\frak{g})\Delta C_0(J,E)$.

(iii) Let $\xi : J\to E$  and $(T(\xi/\frak{g})
)(s,
 \cdot)$,   $(\underline{T}(\frak{g}\xi))(s,\cdot)\in UC (J,E)$) for each $s > 0$. Then $\xi\in  \Delta C_0(J,E)= \Delta C_0(J,X)+i  \Delta C_0(J,X)$.

\

\noindent{\bf{Proof}}. (i) By Lemma 2.5 (ii), $\phi\in C(J,E)$.  We have $(T\phi)(s,\cdot)-(\underline{T}\phi)(s,\cdot)=
v(s,\cdot)\phi_s \in  UC(J,E)$, where $v(s,\cdot)=(\frak{g}(s)\gamma_{2s}-\overline{\frak{g}(s)\gamma_{2s}})= 2i\sin (s^2+2s\cdot)\in BUC(J,
\Cdb)$, $s > 0$. It follows $\phi \sin (s^2+2s(\cdot-s)) \in  UC(J,E)$, $ s > 0$.
 Define $t_k (s)$ by $2s t_k(s) -s^2 = k\pi$, $k\in \Zdb$. Set $K(s,\varepsilon)= \{ t\in J: |\sin (-s^2 +2s t)| < \varepsilon\}= \cup_{k\in \Zdb, t_k(s)\in J}\{ t\in J: |t-t_k(s)| < \delta (\varepsilon)\}$, where $|\sin (2s\delta)|=\varepsilon $. Since $v(s,\cdot-s)\phi \in  UC(J,E)$ and $|v(s,\cdot-s)|\ge \varepsilon$ on $J\setminus K(s,\varepsilon)$, by simple calculations
 \begin{equation}
\label{metric}
 \phi \in  BUC(J\setminus K(s,\varepsilon),E),\qquad s > 0.
 \end{equation}
\noindent  Take $s_1=1$, $s_2=2$.
One can choose  $\varepsilon= \sin (4 \delta)  $ with $ \delta >0 $ small enough such that  $ K(1,\varepsilon)\subset J\setminus K(2,\varepsilon)$. By (16) with $s=1$ and $s=2$, $\phi \in  BUC(J\setminus K(1,\varepsilon),E)$ and $\phi \in  BUC(J\setminus K(2,\varepsilon),E)$.  Since $J \subset (J\setminus K(1,\varepsilon)) \cup (J\setminus K(2,\varepsilon)) $, $\phi\in  BUC(J,E)$.

(ii) If $(T\phi)(s,
 \cdot)\in C_0 (J,E)$, $s >0$, then   $\xi= \frak{g}\phi $  satisfies $\Delta_s\xi \in  C_0(J,X)$  for each $s > 0$. By Lemma 2.5 (i) and (2), $\xi\in \Delta C_0(J,E)$ and  $\phi\in (1/\frak{g})\Delta C_0(J,E)$.

(iii) The assumptions imply $(1/\frak{g})\Delta_s\xi$, $ \frak{g}\Delta_s\xi \in  UC (J,E)$, $s > 0$ by (6). By Lemma 2.1 (iv), $\Delta_s \xi \in C_0 (J,E)$, $s > 0$. By Lemma 2.5 (i) and (2),  $\xi\in  \Delta C_0(J,E)$. $\square$

\

\noindent{\bf {Theorem} 3.2}. Let  $X$, $E$ be  defined by (10).

 (i) If $\phi \in  \Rdb^J$
and
 $(T\phi)(s,\cdot)$  $\in UC(J,\Cdb)$ for each $s> 0$, then $\phi\in BUC(J,\Rdb)$.

(ii)  If $\phi \in  X^J$
and
 $(T\phi)(s,\cdot)$   $\in UC(J,E)$ for each $s> 0$, then $\phi\in BUC(J,X)$.

(iii)  If $f \in  X^J$, $\phi =\frac{f}{\frak{g}}$
with
 $(T
 \phi)(s,\cdot)$  $\in UC(J,E)$ for each $s> 0$, then $f\in \Delta C_0(J,X)$.

\

\noindent{\bf{Proof}}.
 (i) We note that if  $\phi \in  \Rdb^J$, then $\overline{(T\phi)(s,t)}=(\underline{T}\phi)(s,t)$, $s, t \in J $. Since
$(T\phi)(s,\cdot)\in  UC(J,\Cdb)$ for each $s >0$, we conclude  $(\underline{T}\phi)(s,t)=(\overline {T\phi)(s,\cdot)}
\in UC(J,\Cdb)$ for each $s >0$. By Theorem 3.1 (i), $\phi\in BUC(J,\Rdb)$.

(ii) Since $(T\phi)(s,t)=\cos (s^2+2st)\phi (s+t)-\phi (t) + i \sin (s^2+2st)\phi (s+t)$, by (11)(iii) we conclude $\cos (s^2+2s\cdot)\phi_s -\phi,\,   \sin (s^2+2s\cdot)\phi_s \in UC(J,X)$  for each $s >0$. It follows  $(\underline{T}\phi)(s,\cdot)=\cos (s^2+2s\cdot)\phi_s -\phi- i (\sin (s^2+2s\cdot)\phi_s)\in UC(J,E
)$  for each $s >0$. By Theorem 3.1 (i), $\phi\in BUC(J,X)$.

 (iii) By (6), we have $(T\phi)(s,\cdot)= \frac{\Delta_s f}{\frak{g}}= (\Delta_s f\cdot  \cos \cdot^2 -i  \Delta_s f\cdot \sin \cdot^2 ) \in UC(J,E)$ for each $s >0$. By (11)(iii
 ),  $\Delta_s f\cdot  \cos \cdot^2,  \Delta_s f\cdot \sin \cdot^2  \in UC(J,X)$. This implies  $\frak{g} \Delta_s f \in UC(J,E)$, $s > 0$. By Theorem 3.1 (iii), $f\in \Delta C_0(J,X)$.
 $\square$

\

Let   $X$, $E$ be  as defined by (10) and $T$ be as defined by (6), $s > 0$. Denote by

\begin{equation}
\label{metric}
 G=G (J,E)=\{\phi\in E^J: (T\phi) (s,\cdot)\in UC(J,E),\, s > 0\},
\end{equation}
\begin{multline}
\label{metric}
   SC(J,X)= C(J,X)+ (1/\frak{g}) C(J,X),\, \,\,  (TSC(J,X))(s,\cdot)= \{ (T\phi) (s,\cdot): \phi\in SC(J,X)\},\\
(TC(J,X))(s,\cdot)= \{ (T\phi) (s,\cdot): \phi\in C(J,X)\},\,\,\, \Delta_s C(J,X)=\{\Delta_s \phi:  \phi\in C(J,X)\}.\quad
\end{multline}

\noindent Obviously   $G$ is a linear subspace of $E^J$ and by
 Lemma 2.5 (ii), $G \subset C(J,E)$. Also, $ SC(J,X)$ is a linear subspace of $ C(J,E)$.

\

\noindent\textbf{Lemma 3.3}.  Let
  $ SC(J,X)$,
  $ (TSC(J,X))(s,\cdot)$, $\Delta_s C(J,X)$ be  defined by (18).
  \begin{multline}  \label{metric}
\text{\,\,\, (i)}\qquad \,\, \qquad  SC(J,X)=
C (J,X)\oplus (1/\frak{g})C (J,X),\qquad \quad\qquad  \qquad\\
 \text{\, \,(ii)}\qquad \qquad\quad  SC(J,X) \cap UC(J,E)=UC(J,X)\oplus (1/\frak{g})C_0(J,X)\qquad \quad\qquad  \qquad
\end{multline}

(iii)  $(TSC(J,X))(s,\cdot)= (TC(J,X))(s,\cdot)\oplus (1/\frak{g}) \Delta_s C(J,X),\,\,\,  s > 0$.

(iv) If  $\phi= (u +(1/\frak{g})\xi) \in  SC(J,X)$, then $(T\phi)(s,\cdot)\in UC(J,E)$, $s > 0$ if and only if $u \in BUC(J,X) $ and $\xi \in \Delta C_0(J,X) $.

\

 \noindent{\bf{Proof}}. (19)(i): Let $\phi \in  SC(J,X)$, then  $\phi=  u+ (1/\frak{g})\xi$ for  some $ u,\xi\in C(J,X)$. We show that $ u,\xi$ are unique. Indeed,
 by (11) (iv), $\phi = f+i h$ for unique $f, h\in C(J,X)$. So   $\phi(t) = u (t) +(1/\frak{g}(t))\xi (t)= u(t)+ \xi(t) \cos t^2-i\xi (t)\sin t^2 =  f(t)+i h(t)$.   By (11)(i), (ii),   $u(t)+ \xi(t) \cos t^2 =f(t)$, $\xi (t)\sin t^2 =- h(t)$. It follows $\xi (t) = - h(t)/\sin t^2$ for all $t\in J\setminus (t_n)$,  where $t_n = (|n|\pi)^{1/2}$, $n\in \Zdb$. Since $\xi\in C(J,X)$,  $\xi (t_n) = \lim_{t\to t_n} \xi (t)=- \lim_{t\to t_n} h(t)/\sin t^2$ for each $t_n \in J$. This implies that $\xi$ is uniquely defined on $J$. Hence $u (t)=  f(t)-\xi(t) \cos t^2$ is uniquely defined on $ J$. This proves $ u,\xi$ are  unique. Obviously, if $\phi=0$, then $f=h=0$. So,  $\xi (t)\sin t^2 =- h(t)=0$. Hence  $\xi  =0$ and $u= f-\xi \cos \cdot^2=0$.

Now, let $\phi_n= (u_n +(1/\frak{g})\xi_n), \phi= (u +(1/\frak{g})\xi) \in  SC(J,X)$ and let  $\phi_n\to  \phi$ as $n\to \infty$ (in the topology defined by (3)). We show that $u_n\to  u$ and $\xi_n\to  \xi$ as $n\to \infty$. Indeed, let $\phi_n=f_n+ ih_n$, $\phi=f+ ih$ with $f_n, h_n, f,h\in C(J,X)$, $n\in \Ndb$. By (11), $f_n\to f$, $h_n\to h$ as $n\to \infty$. Hence  $\xi_n (t) =- h_n(t)/\sin t^2\to - h(t)/\sin t^2= \xi(t)$ as $n\to\infty$ for each  $t\in J\setminus (t_n)$. This implies   $\xi_n (t) \to\xi(t)$ as $n\to\infty$ for each $t\in J$. So, $u_n(t)=- \xi_n(t) \cos t^2 +f_n(t)\to - \xi(t) \cos t^2 +f(t)=u(t)$. This  proves  (19)(i).

(19) (ii): Let $\psi= (v + (1/\frak{g})\eta)\in  SC(J,X)\cap UC(J,E)$. Then $\psi_s= (v_s + (1/\frak{g}_s)\eta_s)\in  SC(J,X)\cap UC(J,E)$, $s >0$.  By (11)(iii), $v(t+s)+ \eta(t+s)\cos (t+s)^2 $,  $ \eta (t+s)\sin (t+s)^2 $  are uniformly continuous on $J$, $s >0$. It follows  $\eta\in C_0(J,X)$ and  $v\in UC(J,X)$. This and (19)(i) prove (19)(ii).

(iii) Let $\phi= (u + (1/\frak{g})\xi)\in  SC(J,X)$. Since $T$ is linear, $(T\phi)(s,\cdot)=(Tu)(s,\cdot)+ (T(1/\frak{g})\xi)(s,\cdot)=(Tu)(s,\cdot)+ (1/\frak{g})\Delta_s\xi$. If   $(T\phi)(s,\cdot)=0$, $s >0$, then by (6), $\phi = x (1/\frak{g})$ for some $x\in X$. This implies  $u=0$ and $\xi=
 x$. Hence $(Tu)(s,\cdot)=0$, $(1/\frak{g})\Delta_s\xi=0$, $s > 0$.

Now, Let $\phi_n= (u_n +(1/\frak{g})\xi_n), \phi= (u +(1/\frak{g})\xi) \in  SC(J,X)$ and let  $(T\phi_n)(s,\cdot)\to  (T\phi)(s,\cdot)$ as $n\to \infty$ (in the topology defined by (3)), $s > 0$. Noting that $T: C(J,X)\to (TC(J,X))(s,\cdot) $, $s\in J$ is bijective and $(T(1/\frak{g})\eta)(s,\cdot)=(T(1/\frak{g})\tilde{\eta})(s,\cdot)$ if and only if $\tilde{\eta}-\eta=x$ for some $x\in X$, we can replace $(\xi_n), \xi$ by $(\tilde{\xi_n}), \tilde{\xi}$ so that $\tilde{\phi}_n= (u_n +(1/\frak{g})\tilde{\xi_n}), \tilde{\phi}= (u +(1/\frak{g})\tilde{\xi})$ satisfy
$(T\tilde{\phi_n})(s,\cdot)=(T\phi_n)(s,\cdot), (T\tilde{\phi})(s,\cdot)=(T\phi)(s,\cdot)$ and  $\tilde{\phi_n}\to \tilde{ \phi}$ as $n\to \infty$. By (19)(i), $u_n \to u$  and $\tilde {\xi_n}\to \tilde {\xi}$ as $n\to \infty$ implying $(Tu_n)(s,\cdot)\to  (Tu)(s,\cdot)$ and
$ (1/\frak{g}) {\Delta_s \xi_n}=(T(1/\frak{g}) {\xi_n})(s,\cdot)=(T(1/\frak{g})\tilde {\xi_n})(s,\cdot)\to  (T(1/\frak{g})\tilde {\xi})(s,\cdot)=(T(1/\frak{g}){\xi})(s,\cdot)= (1/\frak{g})\Delta_s {\xi}$ as $n\to \infty$. This proves  part (iii).

(iv) By part (iii), $(T\phi)(s,\cdot)= (Tu)(s,\cdot)+\frac{\Delta_s \xi}{\frak{g}}$, $s > 0$ with unique $u, \Delta_s \xi \in C(J,X)$, $s> 0$.
   If $(T\phi)(s,\cdot)\in UC(J,X)$, $s > 0$  but  $(Tu)(s,\cdot)$ or $\frac{\Delta_s \xi}{\frak{g}} \not \in UC(J,E)$ for some  $s > 0$,  then by part (iii) and (8?), we get a contradiction. This proves $(Tu)(s,\cdot)$ and  $\frac{\Delta_s \xi}{\frak{g}} \in UC(J,E)$, $s > 0$. By Lemma 2.1 (ii), $(1/\frak{g})\Delta_s \xi\in UC(J,X)$, $s > 0$ if and only if $\Delta_s \xi\in C_0(J,X)$, $s > 0$. By (2), $\xi\in \Delta C_0 (J,X)$ and by Theorem 3.2 (ii), $u\in BUC(J,X)$.
$\square$

\

 \noindent\textbf{Lemma 3.4}. Let  $X$, $E$ be   defined by (10). If  $\phi \in E^J$  and  $\Delta_s (\frak{g} \phi) \in\frak{g}  UC (J,E)$ for each $ s >0$, then  $ \phi= u+ (1/\frak{g})\xi$  with  $u\in BUC(J,E)$  and $\xi\in \Delta C_0 (J,E)$.

 \

\noindent{\bf{Proof}}.  By Lemma 2.5 (ii), $\phi \in C (J,E)$. Assume first the case $\phi \in SC (J,X)$. By Lemma 3.3 (iv),  $\phi= u+ (1/\frak{g})\xi$ with unique $u \in BUC  (J,X)$, $\xi\in \Delta C_0 (J,X)$.

 In the general case,  let $\phi \in  C (J,E) $. Let  $\rho$ be a distance compatible with
 the topology defined by (3).  We prove that  $\phi = u+ (1/\frak{g})\xi$ for  unique $u, \xi \in C (J,E)$  satisfying
\begin{multline}
\label{metric}
 \rho (u, (1/\frak{g}) C (J,X))=\inf _{y\in  C (J,X)} \rho (u, (1/\frak{g}) y)= \rho (u,0)
  =\rho (u, i (1/\frak{g}) C (J,X)),\qquad\\
  \rho ((1/\frak{g})\xi,  C (J,X))=\inf _{y\in  C (J,X)} \rho((1/\frak{g})\xi,  y)= \rho ((1/\frak{g})\xi,0)=
  \rho ((1/\frak{g})
  \xi, i C (J,X)).\qquad
\end{multline}
\noindent Indeed, assuming $\phi = u+ (1/\frak{g})\xi=\tilde {u}+ (1/\frak{g})\tilde {\xi}$ with
$u,\xi$ and $\tilde {u}, \tilde {\xi}$ satisfying (20), then $u= \tilde {u}+(1/\frak{g})(\tilde {\xi}-\xi) $. If $\tilde {\xi}-\xi\not = 0 $, then $u$ does not satisfy (20).
 This is a contradiction which proves the  uniqueness of $u, \xi$.
  By (11), $u=u_1 + iu_2$ and $\xi=\xi_1 + i\xi_2$ with unique $u_1, u_2, \xi_1,\xi_2\in C(J,X)$. It follows that the map

\noindent $\chi: C(J,E)\to (C(J,X),
C(J,X))+ (1/\frak{g})(C(J,X),C(J,X)), =SC(J,X\times X) $ defined by
 $\chi(\phi) =  (u_1, u_2) + (1/\frak{g})(\xi_1,\xi_2)$ is well defined.
 We have

 $(T\chi(\phi))(s,t)= ((Tu_1)(s,t), (Tu_2)(s,t)) + (1/\frak{g})(\Delta_s
 \xi_1,\Delta_s\xi_2)$.

\noindent Since $ (C(J,X),C(J,X))$  is real linear topological space, as in the proof of (19),

 $ SC(J,X\times X)= (C(J,X),C(J,X))\oplus  (1/\frak{g})(C(J,X),C(J,X))$.

 \noindent As in Lemma 3.3 (iv), if $\psi= (v_1,v_2)+ (1/\frak{g})(\eta_1,\eta_2)\in SC(J,X\times X) $,  $(T(\psi))(s,\cdot)$ is uniformly continuous, $s > 0$ if and only if
$v_1,v_2 \in BUC(J,X) $ and $\eta_1,\eta_2\in  \Delta C_0 (J,X)$. Obviously, if  $v_1$ or $v_2\not \in BUC(J,X) $ or $\eta_1$ or $\eta_2\not\in  \Delta C_0 (J,X)$, then $(T\psi)(s,\cdot)$ is not uniformly continuous. So, if $\psi= \chi(\phi)$, we get   $(T\phi) (s,\cdot)$ is not uniformly continuous and hence $\phi \not \in G$, where  $G$ is defined by (17). It follows $\chi (G)\subset (BUC(J,X),
BUC(J,X))+ (1/\frak{g})(\Delta C_0(J,X),\Delta C_0(J,X)) $. Obviously,
 $ (BUC(J,X),
BUC(J,X))+ (1/\frak{g})(\Delta C_0(J,X),\Delta C_0(J,X)) \subset \chi(G)$ and hence $\chi (G)= (BUC(J,X),
BUC(J,X))+ (1/\frak{g})(\Delta C_0(J,X),\Delta C_0(J,X)) $. This proves that
 $\phi= (u+ (1/\frak{g})\xi) \in G$
 implies $u=u_1+i u_2 \in BUC(J,E) $, $\xi=\xi_1+i\xi_2\in  \Delta C_0 (J,E)$.
$\square$

\

\noindent{\bf{Theorem 3.5}}.  Let  $X$, $E= X\oplus i X$ be  defined by (10) and let $G$ be defined by (17).

$G=BUC (J,E) + \frac{\Delta C_0 (J,E)}{\frak{g}}$
 $ = \{\phi\in E^J
: (T\phi) (s,\cdot)\in BUC(J,E), s > 0\}$.

\

\noindent{\bf{Proof}}. By (6), $\Delta_s (\frak{g}\phi)\in \frak{g} UC(J,E) $ for each $ s > 0$. By Lemma 3.4, $\phi= u+ (1/\frak{g})\xi $ with $u\in BUC(J,E)$ and  $\xi \in \Delta C_0 (J,E
)$. This proves the first equality for $G$.
  The second  equality for $G$   follows by the obvious inclusions

   $ G= BUC (J,E) + \frac{\Delta C_0 (J,E)}{\frak{g}}\subset  \{\phi\in E^J
: (T\phi) (s,\cdot)\in BUC(J,E), s > 0\}\subset G $.
  $\square$

\

\noindent{\bf{Lemma  3.6}}. Let  $X$, $E$ be  defined by (10).

(i)$ C_0(J,E)$ has (S$\Delta$P): If  $u\in E^J$ and $\Delta_s u \in C_0(J,E)$ for each $s > 0$, then $u\in C(J,E)$ and $(u-M_h u)\in C_0(J,E)$ for $h > 0$.

(ii) $\Delta C_0(J,E)$ has $(\Delta)$.

\

\noindent{\bf{Proof}}. (i) By Lemma 2.5 (i),  $u\in UC(J,E)$.
 It follows $\psi: [0,h]\to C_0 (J,E)$, where $\psi (s)= u_s-u$, $s \in [0,h]$ is a continuous function for each $h >0$. Hence the Bochner-Riemann integral  $u-M_h u= -(1/h)\int_0^h (u_s-u)\, ds =  -(1/h)\int_0^h \psi (s)\, ds\in C_0(J,E)$, $h >0$.

(ii) If $f \in L^1_{loc} (J,E)$  with $F= \Delta_s f \in \Delta C_0 (J,E)$ for $s > 0$, then $\Delta_h F \in C_0 (J,E)$ for $h > 0$.
Since  $C_0 (J,E)$ has $(\Delta)$, $\Delta_s (f- M_hf)= (F-M_h F)\in C_0 (J,E)$, $s >0$, $h>0$. By (2),  $(f- M_hf)\in \Delta C_0(J,E)$. This means $\Delta C_0(J,E)$ has $(\Delta)$.
 $\square$

\

\noindent{\bf{Theorem 3.7}}. Let  $X$, $E$ be  defined by (10). Let $f \in E^J$  and $ \Delta_s f \in  \frak{g} UC (J,E)$, $s >0$.

(i) $f\in \frak{g}BUC(J,E) +\Delta C_0 (J,E)$ and
 $\Delta_s f \in\frak{g}\cdot BUC(J,E)$ for each $s>0$.

 (ii) $M_hf \in\Delta  C_0 (J,E)$,  $h > 0$.

(iii) Each of the spaces $\frak{g} BUC(J,E)$,  $\frak{g} UC(J,E)$  has $(\Delta)$.

\

\noindent{\bf{Proof}}. (i)  By Lemma 3.4, $f = \frak{g}u+ v$, where $u\in BUC(J,E)$ and $v\in \Delta C_0(J,E)$. This gives (i).

 (ii)  By part (i), $\Delta_s f= (\Delta_s \frak{g}u+\Delta_sv) \in\frak{g}\cdot BUC(J,E)$ for each $s>0$. So,
$ \Delta_s M_h f=M_h \Delta_s f \in C_0 (J,X)$ by Lemma 2.3 (i) for
each $s>0$, $h>0$. By (2) we get $M_hf \in\Delta  C_0 (J,X)$,  $h > 0$.

(iii) By (i),  $f = \frak{g}u+ v$, where $u\in BUC(J,X)$ and $v\in \Delta C_0(J,X)$. By Lemma 2.3 (i), $M_h \frak{g}u\in  C_0(J,X)$ and by Lemma 3.6 (i), $(v-M_h v)\in  C_0(J,X)$, $h >0$. It follows
$f-M_hf= (\frak{g}u - M_h \frak{g}u+ v-M_h v)\in \frak{g} BUC(J,X)+ C_0(J,X)= \frak{g} BUC(J,X)\subset \frak{g} UC(J,X)$.  This proves $(\Delta)$ for $\frak{g} UC(J,X)$ and $\frak{g} BUC(J,X)$. $\square$

\

\noindent\textbf{4. Difference and $(\Delta)$ properties for  sums.}

\

In this section we study differences and $(\Delta)$ properties for sums of the form $\U+\frak{g}\V$, where $\U, \V\in \{BUC(J,X), UC(J,X)\}$.

\

\noindent\textbf{Proposition 4.1.} Let $ \U, \V\in \{BUC(J,X), UC(J,X)\}$  and $ \W=\U+\frak {g} \V$. Let $f\in L^1_{loc}(J,X)$ and  $\Delta_s f
\in \W$ for each
 $s>0$.
 If $0 \not =\Delta_s f\in \W/C_0$ for each $s > 0$, then

\begin{equation}
\label{metric}
 \Delta_{s}f =  u(s,\cdot)+ \frak{g}
 v(s,\cdot),\text { where\,\,\,} u(s,\cdot)\in \U/C_0,\,\,\,  \frak{g}v(s,\cdot)\in (\frak{g} \V)/C_0
\end{equation}

\noindent
  are  uniquely defined  for each $s  >0$. Moreover

\begin{equation}
\label{metric}
  \text{\, If\,}  \phi  (t)= u(t,0) \text {\, and\,}  F=f-\phi -f(0)\text{\, then\,\,} \Delta_s \phi = u(s,\cdot) \text {\, and \,\,}\Delta_s F= \frak{g} v(s,\cdot).
\end{equation}

\

\noindent\textbf{Proof.}
  By Proposition 2.2, it follows
 $\Delta_{s}f =  u(s,\cdot)+ \frak{g} (t)v(s,\cdot)$ where $u(s,\cdot)\in \U/C_0$,   $\frak{g}v(s,\cdot)\in (\frak{g} \V)/C_0$
are uniquely defined for each $s >0$.  Set  $\phi (t)= u(t,0)$ and   $F= f-\phi- f(0)$. We show that $\Delta_s \phi= u(s,\cdot)$. Indeed,

$\Delta_{s+k} f (t) =  \Delta_k f(t+s)+ \Delta_s f(t)$
$= u (k, t+s)+u (s,t)+\frak{ g}(t+s)v (k, t+s)+ \frak{g}(t) v(s,t)$.

\noindent So, $u(s+k,t)= u(k, t+s)+u(s, t)$  and

$\Delta_s \phi(t)= u(t+s,0)-u (t,0)= u(t,0) + u (s,t)- u(t,0)= u(s,t)$.

\noindent It follows  $\Delta_s \phi \in \U/C_0$, $s >0$ and $\Delta_s F= \Delta_s f-\Delta_s \phi= \frak{g} v(s,\cdot)\in (\frak{g} \V)/C_0$. $\square $

\

\noindent\textbf{Theorem 4.2.} Let $f\in L^1_{loc}(J,X)$ and  $\Delta_s f
\in \U+\frak{g}\cdot \V$ for each
 $s>0$.  Then

 (i) $ f  \in  \Delta \U+\frak{g}\cdot BUC(J,X)$, where $\Delta \U$ is defined by (2).

 (ii) $M_h f \in\Delta \U$

(iii) $\U+\frak{g}\cdot \V$ has ($\Delta$).

\

\noindent\textbf{Proof.} Assume first that $0\not = \Delta_s f \in \W/C_0$  for each $s > 0$. By (21),
   we have   $\Delta_s f=\Delta_s \phi+ \Delta_s F$ with  $ \Delta_s \phi \in \U/ C_0(J,X) $ and $\Delta_s F\in  (\frak{g}\V)/ C_0(J,X)$. By Theorem 3.5 (i), $F = \frak{g} u +v$ for some  $u\in BUC(J,X)$, $v\in \Delta C_0(J,X)$. It follows $\phi = (f-F- f(0))\in L^1_{loc} (J,X)$. Since  $\Delta_s \phi \in  \U/C_0(J,X) $, $s >0$, one gets $\phi \in \Delta \U$ by (2). Hence $f\in \Delta \U+\frak{g} BUC(J,X)$.
The case $\Delta_s f
\in C_0 (J,X)$ for some $s >0$ implies $f\in \Delta C_0 (J,X)$ or $f= \psi+\xi$ with $\psi_s= \psi$ and $\xi\in \Delta C_0 (J,X)$.
Together,  $f\in \Delta \U+ \frak{g}\cdot  BUC(J,X)$. This gives (i).

(ii) We note that if $\psi \in \Delta \U$,  then  $M_h\psi \in \Delta \U$ for each $h > 0$. By part (i) and Theorem 3.7 (ii),  $M_h f=M_h (\phi +F +f(0))\in \Delta \U+ \Delta C_0 (J,X)= \Delta \U$.

(iii) By part (i), if  $f\in L^1_{loc}(J,X)$ and $\Delta_s f
\in \U+\frak{g}\cdot \V$ for each
 $s>0$, then $f=\phi + F+ f(0)$ with $\phi, F$ satisfying (21), (22).   Since $\phi, F\in L^1_{loc} (J,X)$  with  $\Delta_s\phi \in \U$, $\Delta_s F\in  \frak{g}BUC $ and $ \U$ has $(\Delta)$ by [4, Example 3.3], $ \frak{g}BUC $ has $(\Delta)$ by Theorem 3.7 (iii), part (iii) follows. $\square $

\

\end{document}